# Cut-Pursuit Algorithm for Regularizing Nonsmooth Functionals with Graph Total Variation

**Hugo Raguet**
hugo.raguet@gmail.com

**Loïc Landrieu**
IGN, *MATIS*
landrieu.loic@gmail.com

**Abstract**

We present an extension of the cut-pursuit algorithm, introduced by Landrieu and Obozinski (2017), to the *graph total-variation* regularization of functions with a separable nondifferentiable part. We propose a modified algorithmic scheme as well as adapted proofs of convergence. We also present a heuristic approach for handling the cases in which the values associated to each vertex of the graph are multidimensional. The performance of our algorithm, which we demonstrate on difficult, ill-conditioned large-scale inverse and learning problems, is such that it may in practice extend the scope of application of the total-variation regularization.

## 1 Introduction

Landrieu and Obozinski (2017) recently presented a working-set strategy for minimizing differentiable functions regularized by a *total-variation* seminorm structured on a weighted graph. We propose a modified scheme extending the scope of the algorithm to functions with a nondifferentiable part which is separable along the vertices of the graph. Given a finite graph $G \stackrel{\text{def}}{=} (V, E, w)$ with edge weights $w \in \mathbb{R}_+^E$, the problem is to

$$\text{minimize} \quad F: x \mapsto f(x) + \sum_{v \in V} g_v(x_v) + \sum_{(u,v) \in E} w_{(u,v)} |x_u - x_v|, \tag{P1}$$

where $f: \mathbb{R}^V \to \mathbb{R}$ is differentiable, and for all $v \in V$, $g_v: \mathbb{R} \to {]-\infty, +\infty]}$. Our framework allows us to make only weak assumptions on the regularity of the functions $g_v$. In order to handle infinite values and nondifferentiability, we suppose that for all $v \in V$, $g_v$ is *directionally differentiable*; this is detailed later in definition 2.1, using the notions of *domain* and *directional derivative*.

Our algorithm allows us to find *stationary points* of $F$, that is to say points at which all directional derivatives of $F$ are nonnegative. It can be noted already that if all the considered functionals are convex, then the above hypothesis holds, and a stationary point is equivalent to a global minimum.

Our extension of the cut-pursuit algorithm is motivated by the presence of nondifferentiable terms besides the graph total variation in a wide variety of applications. In signal processing for example, the nondifferentiable $\ell_1$ norm fidelity is used for outlier removal, or the denoising of images corrupted by a multiplicative noise, as presented in the works of Nikolova (2004) and Durand et al. (2010), respectively. Wu et al. (2015) make also use of such fidelity terms for the 3D mesh denoising problem, in which noises of multiple natures and scales may coexist.

Alternatively, additional separable nondifferentiable regularization terms can also be considered. On some problems, it is relevant to constrain the values associated to each vertex within a convex subset of the reals, yielding *box constraints*. Another popular regularization is the *fused LASSO*, introduced initially by Tibshirani et al. (2005), where a LASSO-like penalty complements the total-variation regularizer for enforcing solutions which are both sparse and piecewise constant. Gramfort et al. (2013) propose applications in functional magnetic resonance imaging, Becker et al. (2014) in electroencephalography, Omranian et al. (2016) in genetics, and Takayama and Iwasaki (2016) in remote sensing.

In addition, we also consider extending the cut-pursuit algorithm to the cases where the values taken by the vertices belong to a multidimensional space rather than being scalar. This extension is motivated by the success of convex relaxations of the combinatorial Potts model to solve labeling problems, as reported by Nieuwenhuis et al. (2013).

## 1.1  Related Works

Large-scale problems regularized with the graph total variation are typically solved using proximal splitting algorithms; see the general review of Combettes and Pesquet (2008), and the more specific approaches of Couprie et al. (2013) or Raguet and Landrieu (2015). They are however first-order methods, for which convergence is known to be slow, even though preconditioning strategies can help as argued by Pock and Chambolle (2011) and Raguet and Landrieu (2015).

The connection between graph cuts and total variation has been successfully exploited by Chambolle and Darbon (2009) to reformulate the graph total-variation regularization as a *parametric maximum flow* problem when $f$ is a square $\ell_2$ norm. Xin et al. (2016) extend this to the fused LASSO regularization already mentioned, by composing the above method with the proximity operator of the $\ell_1$ norm. They still have to resort to proximal splitting for dealing with more general functionals.

In another line of thought, Bach et al. (2012) have shown that the sparsity of the solution should be exploited computationally to solve large-scale optimization problems faster. *Working-set algorithms* have shown promising results for the convex setting, as demonstrated by Harchaoui et al. (2015).

These ideas are at the heart of the cut-pursuit algorithm originally proposed by Landrieu and Obozinski (2017), which we detail in the following; we refer to their article for a more in-depth discussion on its connection with other works.

## 1.2  Cut-Pursuit for Differentiable Functions

The algorithmic structure of cut-pursuit is summarized in algorithm 1 below. We expose the general principles behind it, recalling for now the situation where there is no nondifferentiable functional besides the graph total variation, reducing to problem P1 where for all $v \in V$, $g_v \stackrel{\text{def}}{=} 0$.

Since the total-variation penalty has a spatially regularizing effect, solutions are expected to be coarse, that is to say they can be expressed as a vector which is piecewise constant with respect to a partition $\mathcal{V}$ of $V$ into few connected components. A key concept of the cut-pursuit algorithm is the *reduced problem*, which is problem P1 constrained on the space of piecewise constant vectors with respect to $\mathcal{V}$; this can be formulated as

$$\text{minimize} \quad F^{(\mathcal{V})} \colon \begin{array}{rcl} \mathbb{R}^{\mathcal{V}} & \longrightarrow & \mathbb{R} \,, \\ \xi & \longmapsto & F(\sum_{U \in \mathcal{V}} \xi_U 1_U) \,, \end{array} \quad (\text{P2})$$

where for all subset $U \subseteq V$, $1_U \in \mathbb{R}^V$ denotes the vector such that for all $v \in V$, $(1_U)_v \stackrel{\text{def}}{=} 1$ if $v \in U$, 0 otherwise. By factorizing finite differences, the graph total-variation term in $F^{(\mathcal{V})}$ becomes $\sum_{(U,U') \in \mathcal{E}} \omega_{(U,U')} |\xi_U - \xi_{U'}|$, where we note the set of adjacent components

$$\mathcal{E} \stackrel{\text{def}}{=} \left\{ (U, U') \in \mathcal{V}^2 \,\middle|\, (U \times U') \cap E \neq \emptyset \right\} \,,$$

and for all $(U, U') \in \mathcal{E}$, the *edge weight* $\omega_{(U,U')} \stackrel{\text{def}}{=} \sum_{(u,v) \in (U \times U') \cap E} w_{(u,v)}$. Since the term $\xi \mapsto f(\sum_{U \in \mathcal{V}} \xi_U 1_U)$ is still differentiable, the reduced problem P2 is structured just as the original problem P1, but over the *reduced graph* $\mathcal{G} \stackrel{\text{def}}{=} (\mathcal{V}, \mathcal{E}, \omega)$, whose vertices are the components in $\mathcal{V}$. Consequently, it should be much easier to solve.

---

**Algorithm 1:** Principle of cut-pursuit; $D \subset \mathbb{R}^V$ is a set of directions adapted to the problem.

**Initialize:** $\mathcal{V} \leftarrow \{V\}$;
**repeat**
    find $\xi^{(\mathcal{V})} \in \mathbb{R}^{\mathcal{V}}$, stationary point of $F^{(\mathcal{V})} \colon \xi \mapsto F(\sum_{U \in \mathcal{V}} \xi_U 1_U)$;
    $x \leftarrow \sum_{U \in \mathcal{V}} \xi_U^{(\mathcal{V})} 1_U$;
    find $d^{(x)} \in D$, minimizing $d \mapsto F'(x, d)$;
    $\mathcal{V} \leftarrow \bigcup_{U \in \mathcal{V}} \{\text{maximal constant connected components of } (d_u^{(x)})_{u \in U}\}$;
**until** $F'(x, d^{(x)}) \geq 0$;
**return** $x$.

---

The cut-pursuit algorithm iteratively refines the partition $\mathcal{V}$, initialized at $\{V\}$. At each iteration, the reduced problem corresponding to the current partition $\mathcal{V}$ is solved, and its solution is used in turn to refine the components of $\mathcal{V}$.

The rationale of the refinement step stems from the structure of the directional derivative $F'(x, d)$ of $F$ at point $x \in \mathbb{R}^V$ in direction $d \in \mathbb{R}^V$. Some calculus shows that $F'(x, d)$ can be expressed as

$$F'(x, d) = \sum_{v \in V} \delta_v(x) d_v + \sum_{(u,v) \in E_{=}^{(x)}} w_{(u,v)} |d_u - d_v| \,, \quad (1)$$

where

$$\delta_v(x) \stackrel{\text{def}}{=} \nabla_v f(x) + \sum_{\substack{(e,u) \in E \times V \\ e = (u,v) \text{ or } (v,u)}} w_e \, \text{sign}(x_v - x_u) \,,$$

sign: $\mathbb{R} \mapsto \{-1, 0, +1\} \colon t \mapsto -1$ if $t < 0$, $0$ if $t = 0$ and $+1$ if $t > 0$, and $E_{=}^{(x)} \stackrel{\text{def}}{=} \{(u, v) \in E \,|\, x_u = x_v\}$ is the set of edges whose vertices share the same value. The first sum in equation 1 consists in unary terms, in which the sign of $\delta_v(x)$ determines whether the value of each vertex should tend to decrease or increase. The second sum consists in binary terms, encouraging the values at neighboring vertices to evolve in unison.

The refinement step does not require finding a "steepest descent" direction, but merely refining the current partition into a new one, thus adding relevant degrees of liberty to the next reduced problem. The goal is to split the current components into groups of vertices tending to increase together or decrease together, while taking coupling terms into account. Such a split can thus be encoded as a direction in the set $\{-1, +1\}^V$. We thus look for a *steepest binary direction*

$$\text{find} \quad d^{(x)} \in \argmin_{d \in \{-1, +1\}^V} F'(x, d) \,, \tag{P3}$$

which can be solved by finding the minimum cut in an appropriate flow graph. The refined partition $\mathcal{V}$ is then defined by splitting each component $U$ of the current partition according to the constant connected components of $\left(d_u^{(x)}\right)_{u \in U}$.

Beyond the computational efficiency of the cut-pursuit algorithm, the main result of Landrieu and Obozinski (2017) is the optimality certificate, which states that if $x$ is a solution of a reduced problem P2, and that the steepest binary direction problem P3 induces no refinement of the current partition then $x$ is a solution of the main problem P1. In consequence, the algorithm converges in a finite number of steps to such a solution. In practice, since the number of components of $\mathcal{V}$ increases rapidly and the final partition is expected to be coarse, only a few iterations are needed.

### 1.3  Contributions

If we assume now that $F$ has a nondifferentiable part other than the graph total variation, the analysis above does not stand because one cannot decompose the directional derivative into unary and binary contributions as in equation 1.

Assuming that such a nondifferentiable part is separable along vertices, a similar decomposition can be achieved, but the unary terms $\delta_v$ would depend on the sign of $d_v$. In contrast to what happens with the differentiable term, where each vertex either tends to increase or decrease, it is now possible that both directions $+1$ and $-1$ are unfavorable for some vertices. Thus, if one wants to keep the principle of the cut-pursuit algorithm for the regularization of nondifferentiable functions, it seems necessary to search for descent directions within the set $\{-1, 0, +1\}^V$ when refining the partition.

In this paper, we provide a simple theoretical framework allowing to deal with directional derivatives in possibly noncontinuous settings. Then, we show that refining partitions with descent directions within $\{-1, 0, +1\}^V$ is actually sufficient in order to retain the optimality certificate of cut-pursuit with separable nondifferentiable terms. Since the resulting problem can also be solved via a minimum cut in an adapted flow graph, this finally enables the use of the powerful cut-pursuit approach on a large class of problems as introduced above.

Considering a problem in which each vertex takes multidimensional values, unit vectors encoding a descent direction at a vertex are not restricted to the finite $\{-1, +1\}$ set. There is actually an infinity of such unit vectors, and searching for a steepest unit descent direction is intractable. However, we propose some heuristics, and show numerically that by restricting the search to a small set of well-chosen directions, one can still apply the cut-pursuit approach, drastically outperforming traditional proximal schemes.

## 2  Extending Cut-Pursuit

Since the cut-pursuit relies on directional derivatives, we start with some definitions allowing us to manipulate them with the necessary degree of generality.

**Definition 2.1.** Let $\Omega$ be a real vector space, and $h\colon \Omega \to\,]-\infty,+\infty]$. The *domain* of $h$ is $\operatorname{dom} h \stackrel{\text{def}}{=} \{x \in \Omega \mid h(x) < +\infty\}$. Given $x \in \operatorname{dom} h$ and $d \in \Omega$, we say that $h$ admits a *directional derivative at point $x$ in direction $d$* if the quantity $h'(x, d) \stackrel{\text{def}}{=} \lim_{t \downarrow 0} \frac{h(x+td) - h(x)}{t}$ exists in $]-\infty, +\infty]$. Finally, we say that $h$ is *directionally differentiable* if it admits a directional derivative at every point of its domain and in every direction.

Our definition of directional derivatives would be standard, if it were not for infinite values. It is easy to establish that sums of directionally differentiable functions are directionally differentiable, and that the directional derivative of the sum is the sum of the directional derivatives.

Interestingly, in dimension one, that is $\Omega \stackrel{\text{set}}{=} \mathbb{R}$, one can easily show that *directional differentiability implies lower semicontinuity over the domain*. The reciprocal does not hold; consider $h\colon x \mapsto x \sin(1/x)$ and $h(0) \stackrel{\text{def}}{=} 0$, which is continuous at 0 but not directionally differentiable at this point. Note also that this does not hold in larger dimensions, as shown by the counterexample $h\colon \mathbb{R}^2 \to \mathbb{R}\colon (x, y) \mapsto -x^2 y/(x^4 + y^2)$ and $h(0, 0) \stackrel{\text{def}}{=} 0$. Indeed, for all $(d, e) \in \mathbb{R}^2$ and $t \in \mathbb{R}$, $h(t(d, e)) = -td^2 e/(t^2 d^4 + e^2)$, and thus $h'((0, 0), (d, e)) = 0$. However, for all $x \in \mathbb{R}$, $h(x, x^2) = -1/2 < 0$, hence $h$ is not lower semicontinuous at 0.

Finally, a special case of interest is convexity, *automatically ensuring directional differentiability*. In the usual case with only finite values, this is typically shown by Hiriart-Urruty and Lemaréchal (2004, Part D), recalling that directional derivative of $h$ at point $x$ in direction $d$ is nothing but the right derivative at 0 of the unidimensional functional $\ell\colon t \mapsto h(x + td)$, which is also convex, admitting thus the right derivative in question. If now $h$ admits positively infinite values but $x$ belongs to its domain, either $\ell(t) = +\infty$ for all $t > 0$, in which case $h'(x, d) = +\infty$, or $\ell(t) < +\infty$ for some $t > 0$, in which case $\ell$ is finite over $[0, t]$ by convexity, reducing to the above case with only finite values.

Convexity is of particular importance because many applications would use convex optimization algorithms for solving the reduced problem P2. However, we underline that it is not a requirement, and that the cut-pursuit algorithm can be perfectly applied on nonconvex problems, provided that solutions of the reduced problems can be found.

In the remainder of this section, we first describe our method for extending the cut-pursuit algorithm, and the rationale behind it. We then further justify this rationale by providing a convergence proof. Subsequently, we specify some practical implementation details. Finally, we give an efficient heuristic for dealing with a similar setting where the values at each vertex are multidimensional.

## 2.1 Steepest Ternary Direction

As stated in § 1.3, nondifferentiable terms in $F$ prevent convenient decomposition of the directional derivatives as equation 1. However, an important property of the directional derivative is that it is *positively homogeneous*.

**Proposition 2.1.** *Let $\Omega$ be a real vector space, and $h\colon \Omega \to\,]-\infty,+\infty]$. If $h$ admits a directional derivative at $x \in \operatorname{dom} h$ and direction $d \in \Omega$, then for all $\lambda \geq 0$, it also admits a directional derivative at $x$ in direction $\lambda d$, and $h'(x, \lambda d) = \lambda h'(x, d)$.*

*Proof.* If $\lambda = 0$, then directly $h'(x, \lambda d) = 0$; suppose then that $\lambda > 0$. For all $t > 0$, write

$$\frac{h(x + t(\lambda d)) - h(x)}{t} = \lambda \frac{h(x + (t\lambda)d) - h(x)}{t\lambda},$$

which holds even if $h(x + t\lambda d)$ is infinite. Now, as $t$ tends to zero with positive values, so does $t\lambda$, thus the above quantity tends to $\lambda h'(x, d)$. ∎

This positive homogeneity is especially useful in our case where the nondifferentiable terms are separable into unidimensional functionals. Decomposition similar to equation 1 can be achieved, although multiplicative terms $\delta_v$ now depend on the sign of the corresponding direction coordinate.

**Proposition 2.2.** *Under our assumptions, for all $x \in \mathrm{dom}\, F$ and for all $d \in \mathbb{R}^V$, $F$ admits a directional derivative at $x$ in direction $d$, equal to*

$$F'(x, d) = \sum_{\substack{v \in V \\ d_v > 0}} \delta_v^+(x) d_v + \sum_{\substack{v \in V \\ d_v < 0}} \delta_v^-(x) d_v + \sum_{(u,v) \in E_=^{(x)}} w_{(u,v)} |d_u - d_v|, \qquad (2)$$

*where for all $v \in V$, we define*

$$\delta_v^+(x) \stackrel{\mathrm{def}}{=} \nabla_v f(x) + g_v'(x_v, +1) + \sum_{\substack{(e,u) \in E \times V \\ e=(u,v) \text{ or } (v,u)}} w_e \,\mathrm{sign}(x_v - x_u),$$

*and*

$$\delta_v^-(x) \stackrel{\mathrm{def}}{=} \nabla_v f(x) - g_v'(x_v, -1) + \sum_{\substack{(e,u) \in E \times V \\ e=(u,v) \text{ or } (v,u)}} w_e \,\mathrm{sign}(x_v - x_u).$$

*Proof.* The existence of the directional derivative is provided as the sum of finite or positively infinite limits. The directional derivative of $f$ and of $(x_1, x_2) \mapsto |x_1 - x_2|$ at a point where $x_1 \neq x_2$ is given by differentiability. Then, the directional derivative of $(x_1, x_2) \mapsto |x_1 - x_2|$ at a point where $x_1 = x_2$ in the direction $(d_1, d_2)$ is easily found by writing in that case for all $t > 0$, $|x_1 + td_1 - x_2 - td_2|/t = |d_1 - d_2|$. Finally, let $v \in V$. If $d_v = 0$ then $g_v'(x, d_v) = 0$. Otherwise, using proposition 2.1 we get $g_v'(x, d_v) = g_v'(x, |d_v|\,\mathrm{sign}(d_v)) = |d_v| g_v'(x, \mathrm{sign}(d_v))$ that is to say $g_v'(x, +1) d_v$ if $d_v > 0$ and $-g_v'(x, -1) d_v$ if $d_v < 0$. ∎

In contrast to the differentiable case, it is now possible that for some vertices, neither increasing $+1$ nor decreasing $-1$ direction is favorable, when looking for convenient descent directions. In this case, such vertices are inclined not to change their value, that is to say the null direction $0$ should be favored; this leads to the *steepest ternary direction* problem,

$$\text{find} \quad d^{(x)} \in \argmin_{d \in \{-1, 0, +1\}^V} F'(x, d), \qquad (\text{P4})$$

where for all $x \in \mathrm{dom}\, F$ and $d \in \{-1, 0, +1\}^V$,

$$F'(x, d) = \sum_{\substack{v \in V \\ d_v = +1}} \delta_v^+(x) - \sum_{\substack{v \in V \\ d_v = -1}} \delta_v^-(x) + \sum_{(u,v) \in E_=^{(x)}} w_{(u,v)} |d_u - d_v|.$$

Observe that since $F'(x, 0) = 0$, any solution $d^{(x)}$ of problem P4 must satisfy $F'(x, d^{(x)}) \leq 0$. Similarly to its binary counterpart, the steepest ternary direction corresponds to a minimum cut in a suitable flow graph, represented in figure 1, which we note $G_{\mathrm{flow}}^{(x)} = (V_{\mathrm{flow}}, E_{\mathrm{flow}}^{(x)}, c^{(x)})$. The vertex set is $V_{\mathrm{flow}} = (V \times \{1, 2\}) \cup \{s, t\}$, $s$ and $t$ being respectively the specific *source* and *sink* vertices; we also use the convenient notation $v^{(k)}$ for $(v, k) \in V \times \{1, 2\}$. The edge set is defined by

$$E_{\text{flow}}^{(x)} \stackrel{\text{def}}{=} \bigcup_{v \in V} \left\{ (s, v^{(1)}), (v^{(1)}, v^{(2)}), (v^{(2)}, t) \right\} \cup \bigcup_{\substack{(u,v) \in E_{=}^{(x)} \\ k \in \{1,2\}}} \left\{ (u^{(k)}, v^{(k)}), (v^{(k)}, u^{(k)}) \right\}.$$

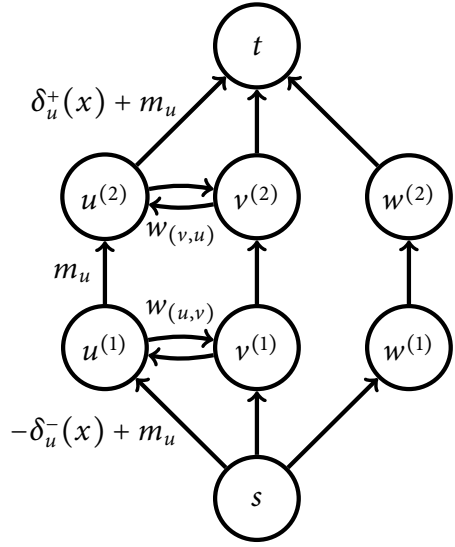

Figure 1: Schematic representation of the flow graph $G_{\text{flow}}^{(x)}$ for the steepest ternary direction problem P4. In this illustration, $x_u = x_v \neq x_w$.

In accordance with figure 1, the edges defined in the left term are called *vertical*, while the edges defined in the right term are called *horizontal*.

The associated capacities $c^{(x)} \in \mathbb{R}_+^{|E_{\text{flow}}^{(x)}|}$ are defined, for the horizontal edges, for all $(u, v) \in E_{=}^{(x)}$ and $k \in \{1, 2\}$, by $c_{(u^{(k)}, v^{(k)})}^{(x)} \stackrel{\text{def}}{=} c_{(v^{(k)}, u^{(k)})}^{(x)} \stackrel{\text{def}}{=} w_{(u,v)}$; and for the vertical edges, for all $v \in V$, by

$$c_{(s, v^{(1)})}^{(x)} \stackrel{\text{def}}{=} -\delta_v^-(x) + m_v, \qquad c_{(v^{(1)}, v^{(2)})}^{(x)} \stackrel{\text{def}}{=} m_v \qquad \text{and} \qquad c_{(v^{(2)}, t)}^{(x)} \stackrel{\text{def}}{=} \delta_v^+(x) + m_v,$$

where $m_v \stackrel{\text{def}}{=} \max(0, \delta_v^-(x), -\delta_v^+(x))$, $\delta_v^-(x)$ and $\delta_v^+(x)$ being defined in proposition 2.2; note that our definition of directional derivatives implies that $\delta_v^-(x) < +\infty$ and $-\delta_v^+(x) < +\infty$. The definition of $m_v$ ensures that all capacities are nonnegative, although potentially infinite. An additional benefit is that for each $v \in V$, at least one of $c_{(s, v^{(1)})}^{(x)}$, $c_{(v^{(1)}, v^{(2)})}^{(x)}$ and $c_{(v^{(2)}, t)}^{(x)}$ is zero, allowing for faster computation of the minimum cut via an augmenting path algorithm, such as the one of Boykov and Kolmogorov (2004).

It can also be noted that this flow graph is similar to the multistage structure proposed by Ishikawa (2003), with one fewer stage and no infinite so-called *constraint edges*; this is once again favorable to augmenting path algorithms.

**Proposition 2.3.** *Problem P4 can be solved by finding a minimum cut in the graph $G_{\text{flow}}^{(x)}$.*

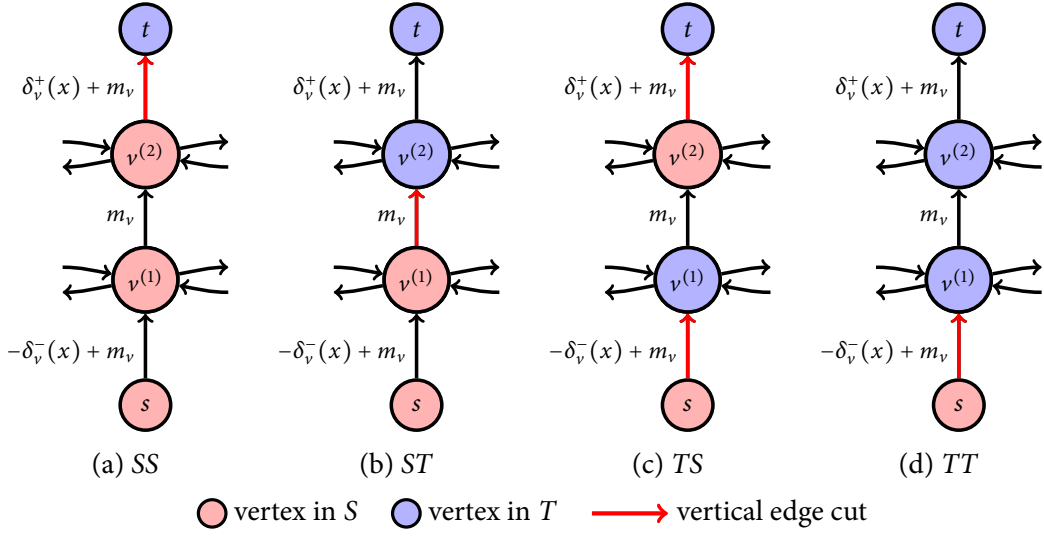

Figure 2: Illustration of the configurations *SS*, *ST*, *TS*, and *TT*.

*Proof.* A cut is an ordered partition $(S, T)$ of $V_{\text{flow}}$ such that $s \in S$ and $t \in T$; we denote the *cut value* by $\text{cut}(S, T) \stackrel{\text{def}}{=} \sum_{(u,v) \in E^{(x)}_{\text{flow}} \cap (S \times T)} c^{(x)}_{(u,v)}$. It is also convenient to define the *nonoriented edge weight* between two subsets $U, U' \subset V$ by $w(U, U') = \sum_{(u,v) \in E^{(x)}_{\leq} \cap ((U \times U') \cup (U' \times U))} w_{(u,v)}$.

Let $(S, T)$ be a cut. Given a vertex $v \in V$, there are four possible configurations depending on which set contains $s$, $v^{(1)}$, $v^{(2)}$ and $t$, as illustrated in figure 2. Accordingly, we define the sets *SS*, *ST*, *TS*, *TT* $\subseteq V$, where for $A, B \in \{S, T\}$, $AB \stackrel{\text{def}}{=} \{v \in V \mid v^{(1)} \in A \text{ and } v^{(2)} \in B\}$.

We argue that if the set *TS* is nonempty, then one can find an alternative cut with a smaller or equal value than $(S, T)$ and for which no vertex in $V$ presents the configuration *TS*. Indeed, for any vertex $v \in TS$, both vertical edges $(s, v^{(1)})$ and $(v^{(2)}, t)$ are cut by $(S, T)$. Now, the number of horizontal edges cut is determined by the configuration of the neighbors: one for each neighbor in *SS* or *TT*, two for each neighbor in *ST*, and zero for each neighbor also in *TS*. We can thus write

$$\text{cut}(S, T) = \sum_{v \in TS} \left( c^{(x)}_{(s,v^{(1)})} + c^{(x)}_{(v^{(2)},t)} \right) + w(TS, SS) + w(TS, TT) + 2\, w(TS, ST) + K,$$

where $K$ is a constant which does not depend on the configuration of the vertices in *TS*. Let us now consider the cut $(S', T')$ in which all vertices of *TS* are switched to configuration *SS*; more precisely $S' = S \cup (TS \times \{1\})$ and $T' = T \smallsetminus (TS \times \{1\})$. Likewise, we consider $(S'', T'')$ for which the vertices of *TS* are switched to configuration *TT*; $S'' = S \smallsetminus (TS \times \{2\})$ and $T' = T \cup (TS \times \{2\})$. We can write the value of these cuts as

$$\text{cut}(S', T') = \sum_{v \in TS} c^{(x)}_{(v^{(2)},t)} + w(TS, ST) + 2\, w(TS, TT) + K,$$
$$\text{cut}(S'', T'') = \sum_{v \in TS} c^{(x)}_{(s,v^{(1)})} + w(TS, ST) + 2\, w(TS, SS) + K.$$

We now show that at least one of $\text{cut}(S', T')$ or $\text{cut}(S'', T'')$ is smaller or equal to $\text{cut}(S, T)$.

Indeed, we have the following:

$$\text{cut}(S, T) - \text{cut}(S', T') = \sum_{v \in TS} c^{(x)}_{(s,v^{(1)})} + w(TS, ST) + w(TS, SS) - w(TS, TT)$$

$$\text{cut}(S, T) - \text{cut}(S'', T'') = \sum_{v \in TS} c^{(x)}_{(v^{(2)},t)} + w(TS, ST) - w(TS, SS) + w(TS, TT).$$

Since all edge weights are positive, at least one of $\text{cut}(S,T) - \text{cut}(S',T')$ or $\text{cut}(S,T) - \text{cut}(S'',T'')$ must be positive, unless simultaneously all vertical edges have zero capacity, $w(TS, ST) = 0$, and $w(TS, SS) = w(TS, TT)$. In which case, all three cuts $(S, T)$, $(S', T')$, and $(S'', T'')$ have the same value.

We have proven that for any cut $(S, T)$ for which $TS$ is nonempty, we can find an alternative cut for which it is empty and with smaller or equal value. Consequently, for a minimum cut $(S, T)$ of $G^{(x)}_{\text{flow}}$, the vertices in $TS$ can be switched to configuration $SS$ without changing the resulting cut value. Finally, to each such cut $(S, T)$ satisfying $TS = \emptyset$, we associate a direction $d \in \{-1, 0, 1\}^V$ defined for all $v \in V$ by $d_v \stackrel{\text{def}}{=} -1$ if $v \in TT$, $0$ if $v \in ST$, and $+1$ if $v \in SS$. This mapping is one-to-one, and satisfies

$$\text{cut}(S, T) = \sum_{v \in SS} c^{(x)}_{(s,v^{(1)})} + \sum_{v \in ST} c^{(x)}_{(v^{(1)},v^{(2)})} + \sum_{v \in TT} c^{(x)}_{(v^{(2)},t)}$$
$$+ w(SS, ST) + w(ST, TT) + 2w(SS, TT),$$
$$= F'(x, d) + \sum_{v \in V} m_v.$$

Since the term $\sum_{v \in V} m_v$ does not depend on $d$, minimizing the value of the cut amounts to minimizing $d \mapsto F'(x, d)$. ∎

## 2.2 Convergence Proof

We now turn to the convergence of algorithm 1 towards a stationary point of $F$. In our context, a strictly negative directional derivative is called a *strict descent direction*, and a point $x \in \text{dom}\, F$ is called *stationary* if it admits no strict descent direction.

The convergence proof relies on the same kind of optimality certificate than the one used in the original cut-pursuit paper by Landrieu and Obozinski (2017), for regularization of differentiable functionals. Indeed, the steepest ternary direction at a point $x$ not only indicates a refinement of the partition $\mathcal{V}$, it also allows us to determine the optimality of $x$ as a solution of the main problem P1, even though it is not a steepest descent direction in general.

However, while the proofs for differentiable functionals consider the flows in $G^{(x)}_{\text{flow}}$ and relate to atomic gauge theory, ours are more elementary. They focus on making explicit the fact that, as long as the current iterate $x$ is not stationary, steepest ternary descent directions induce further refinements of the partition $\mathcal{V}$. To do so, we introduce the concept of *sign-segregation* of a vector defined over the graph $G$, which conveys the idea that neighboring vertices with the same sign must also have the same value; this will be useful for manipulating descent directions.

**Definition 2.2.** We say that a vector $d \in \mathbb{R}^V$ is *sign-segregated* over $G$, if for all $(u, v) \in E$, $\text{sign}(d_u) = \text{sign}(d_v)$ implies $d_u = d_v$.

Consider then the following lemma.

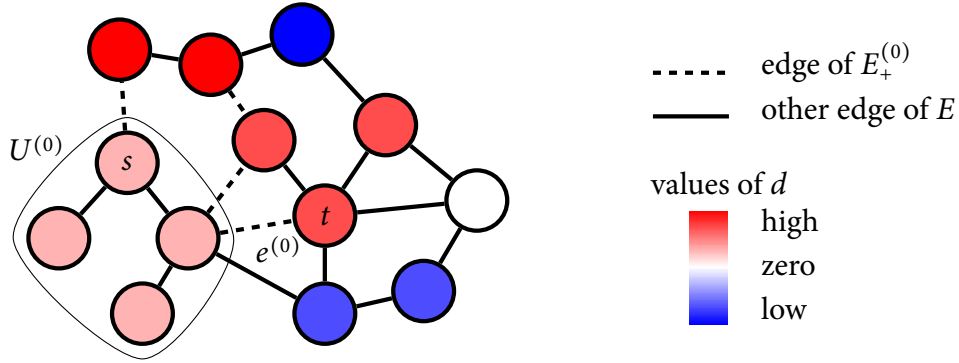

Figure 3: Illustration for the proof of lemma 2.1; $d^{(0)}$ is schematized in colors, positive coordinates in red, negative in blue.

**Lemma 2.1.** *Let $x \in \text{dom } F$. If $F$ admits a strict descent direction at point $x$, then it admits a strict descent direction which is sign-segregated over $G$.*

*Proof.* Let $d^{(0)} \in \mathbb{R}^V$ be a strict descent direction for $F$ at $x \in \mathbb{R}^V$. We construct recursively a strict descent direction which is sign-segregated. Suppose first that the set of edges where $d^{(0)}$ violates the sign-segregation with positive values, $E_+^{(0)} \stackrel{\text{def}}{=} \{(u,v) \in E \,|\, \text{sign}(d_u^{(0)}) = \text{sign}(d_v^{(0)}) = +1 \text{ and } d_u^{(0)} \neq d_v^{(0)}\}$, is nonempty.

The set of values of $d^{(0)}$ at the involved vertices, that is to say $\{d_v^{(0)} \in \mathbb{R} \,|\, \exists u \in V : (u,v) \in E_+^{(0)} \text{ or } (v,u) \in E_+^{(0)}\}$, is then nonempty and has a smallest element $d_s^{(0)}$ associated to vertex $s$. Let then $U^{(0)} \subset V$ be the *largest connected component* in $G$ containing $s$ over which $d^{(0)}$ is constant equal to $d_s^{(0)}$. Now, the set of positive values of $d^{(0)}$ at vertices neighboring $U^{(0)}$, that is to say $\{d_v^{(0)} \in \mathbb{R} \,|\, v \in V \smallsetminus U^{(0)} : d_v^{(0)} > 0 \text{ and } \exists u \in U^{(0)} : (u,v) \in E \text{ or } (v,u) \in E\}$, is nonempty since $s \in U^{(0)}$ is a vertex involved in an edge within $E_+^{(0)}$. Again, it has a smallest element $d_t^{(0)}$ associated to vertex $t \in V \smallsetminus U^{(0)}$ and edge $e^{(0)} \in E$ connecting it to $U^{(0)}$. By construction, it holds that for all $u \in U^{(0)}$, $d_u^{(0)} = d_s^{(0)}$, that $d_t^{(0)} > d_s^{(0)} > 0$ and that $e^{(0)} \in E_+^{(0)}$. We show in the following that we can modify $d^{(0)}$ by substituting the value $d_s^{(0)}$, shared by all the $U^{(0)}$ vertices, either by $d_t^{(0)}$ or by 0, hence decreasing the number of edges on which sign-segregation is violated, while keeping a strict descent direction.

According to proposition 2.2, the expression of the directional derivative $F'(x, d^{(0)})$, in function of the value $d_s^{(0)}$ shared by all the $U^{(0)}$ vertices, is, up to a constant depending only on the values at the $V \smallsetminus U^{(0)}$ vertices,

$$d_s^{(0)} \mapsto \sum_{u \in U^{(0)}} \delta_u^+(x) d_s^{(0)} + \sum_{\substack{(e,u,v) \in E_=^{(x)} \times U^{(0)} \times V \\ e=(u,v) \text{ or } (v,u)}} w_e \, \text{sign}(d_s^{(0)} - d_v^{(0)})(d_s^{(0)} - d_v^{(0)}) \,.$$

This expression, with of the quantities $\delta_u^+(x)$, holds as long as $d_s^{(0)}$ is nonnegative. Let now $(u,v) \in U^{(0)} \times V$ such that $(u,v) \in E_=^{(x)}$ or $(v,u) \in E_=^{(x)}$, and proceed by case analysis. If $v \in U^{(0)}$, then $d_v^{(0)} = d_s^{(0)}$. If $v \in V \smallsetminus U^{(0)}$, then $d_v^{(0)} \neq d_s^{(0)}$, otherwise contradicting the definition of $U^{(0)}$, $d_v^{(0)} < d_s^{(0)} \implies d_v^{(0)} \leq 0$, otherwise contradicting the definition of $d_s^{(0)}$, and $d_v^{(0)} > d_s^{(0)} \implies d_v^{(0)} \geq d_t^{(0)}$, otherwise contradicting the definition of $d_t^{(0)}$. We deduce that the quantity $\mathrm{sign}(d_s^{(0)} - d_v^{(0)})$ remain constant for $d_s^{(0)}$ ranging in $]0, d_t^{(0)}[$; moreover, each term of the sum is continuous at 0 and $d_t^{(0)}$. Altogether, the expression is linear in $d_s^{(0)}$ over the range $[0, d_t^{(0)}]$, with coefficient $a^{(0)} \stackrel{\text{def}}{=} \sum_{u \in U^{(0)}} \delta_u^+(x) + \sum_{\substack{(e,u,v) \in E_=^{(x)} \times U^{(0)} \times V \\ e=(u,v) \text{ or } (v,u)}} w_e \, \mathrm{sign}(d_s^{(0)} - d_v^{(0)})$, and reaches its extrema at 0 and $d_t^{(0)}$.

In accordance, we define $d^{(1)}$ such that for all $v \in V \setminus U^{(0)}$, $d_v^{(1)} \stackrel{\text{def}}{=} d_v^{(0)}$, and for all $v \in U^{(0)}$, $d_v^{(1)} \stackrel{\text{def}}{=} 0$ if $a^{(0)} \geq 0$, $d_t^{(0)}$ otherwise. This ensures that $F'(x, d^{(1)}) \leq F'(x, d^{(0)})$, while $d^{(1)}$ does not violate sign-segregation over the edge $e^{(0)}$, because either $d_s^{(1)} = d_t^{(1)}$, or $\mathrm{sign}(d_s^{(1)}) = 0$ and $\mathrm{sign}(d_t^{(1)}) = +1$.

We can then proceed with $E_+^{(1)} \subseteq E_+^{(0)} \setminus e^{(0)}$; doing so recursively provides a sequence of strict descent directions, violating sign-segregation with positive values over a strictly decreasing number of edges. The same recursion can be applied *mutatis mutandis*, to take care of edges over which sign-segregation is violated with negative values. The total number of edges being finite, the desired sign-segregated strict descent direction is obtained after a finite number of recursions. ■

Sign-segregation facilitates the proof of the following fundamental proposition.

**Proposition 2.4.** *Let $x \in \mathrm{dom}\, F$. If $F$ admits a strict descent direction at point $x$, then it admits a strict descent direction in the set $\{-1, 0, +1\}^V$.*

*Proof.* Let $d \in \mathbb{R}^V$ be a strict descent direction for $F$ at $x \in \mathbb{R}^V$; thanks to [lemma 2.1](), we can assume that $d$ is sign-segregated. There exists then a partition $\mathcal{U}$ of $V$ such that $d$ is constant over each component of $\mathcal{U}$, and takes different signs over two components which are neighbors in $G$.

Now, by splitting the absolute difference of $a, b \in \mathbb{R}$ as $|a - b| = \mathrm{sign}(a - b)a + \mathrm{sign}(b - a)b$, we can rewrite the directional derivative [equation 2]() by regrouping the vertices according to the partition $\mathcal{U}$, yielding

$$F'(x, d) = \sum_{U \in \mathcal{U}} \sum_{u \in U} \delta_u(x, \mathrm{sign}(d_u)) d_u + \sum_{\substack{(e,u,v) \in E_=^{(x)} \times U \times V \\ e=(u,v) \text{ or } (v,u)}} w_e \, \mathrm{sign}(d_u - d_v) d_u \,,$$

where for all $u \in V$ and $s \in \{-1, 0, +1\}$, $\delta_u(x, s) \stackrel{\text{def}}{=} \delta_u^+$ if $s = +1$, $\delta_u^-$ otherwise. If $\zeta$ is the vector of $\mathbb{R}^\mathcal{U}$ such that $d = \sum_{U \in \mathcal{U}} \zeta_U 1_U$, this can be factorized as $F'(x, d) = \sum_{U \in \mathcal{U}} \Delta_U(x, d) \zeta_U$, where $\Delta_U(x, d) \stackrel{\text{def}}{=} \sum_{u \in U} \delta_u(x, \mathrm{sign}(d_u)) + \sum_{\substack{(e,u,v) \in E_=^{(x)} \times U \times V \\ e=(u,v) \text{ or } (v,u)}} w_e \, \mathrm{sign}(d_u - d_v)$.

Since $F'(x, d) < 0$, there exists at least one component $U \in \mathcal{U}$ such that $\Delta_U(x, d) \zeta_U < 0$. Let $U$ be such a component, and let $(u, v) \in U \times V$ such that $(u, v) \in E_=^{(x)}$ or $(v, u) \in E_=^{(x)}$. If $v \in U$, we recall simply that $d_v = d_u = \zeta_U$ and $\mathrm{sign}(d_u - d_v) = 0$. Now if $v \in V \setminus U$, then by sign-segregation we know that $\mathrm{sign}(d_u) \neq \mathrm{sign}(d_v)$, and then it follows by case analysis that $\mathrm{sign}(d_u - d_v) = \mathrm{sign}(d_u - 0) = \mathrm{sign}(\zeta_U)$. We deduce that $\Delta_U(x, \mathrm{sign}(\zeta_U) 1_U) = \Delta_U(x, d)$; calling on the above decomposition, we conclude $F'(x, \mathrm{sign}(\zeta_U) 1_U) = \Delta_U(x, d) \, \mathrm{sign}(\zeta_U) < 0$. ■

**Remark 2.1.** The proof of [proposition 2.4]() actually shows a slightly stronger conclusion, namely that $F$ admits a strict descent direction in the set $\{-1, 0\}^V \cup \{0, +1\}^V \subset \{-1, 0, +1\}^V$.

In order to use the above results for proving the termination and correctness of [algorithm 1](), one should ensure that a stationary point of each reduced problem exists and can be found.

In particular, note that when initializing the algorithm with $\mathcal{V} \stackrel{\text{set}}{=} \{V\}$ as suggested, the mere existence of feasible points for the first reduced problem requires that $\cap_{v \in V} \operatorname{dom} g_v \neq \emptyset$. Of course, it is still possible to initialize with a finer partition if necessary, but these considerations are all problem-dependent. For the scope of the present article, we thus assume the necessary existence properties.

**Corollary 2.1.** *Under our assumptions, algorithm 1 with $D \stackrel{\text{set}}{=} \{-1, 0, +1\}^V$ finds a stationary point of the main problem P1.*

*Proof.* We claim that, at the last step of each iteration, if the iterate $x = \sum_{U \in \mathcal{V}} \xi_U^{(\mathcal{V})} 1_U$ is not a stationary point, then the cardinal of $\mathcal{V}$ is strictly increased. Indeed, supposing otherwise means that for all $U \in \mathcal{V}$, the set of maximal constant connected components of $(d_u^{(x)})_{u \in U}$ is reduced to $\{U\}$, or in other words, $d^{(x)}$ is constant over $U$; thus there exists $\zeta \in \mathbb{R}^{\mathcal{V}}$ such that $d^{(x)} = \sum_U \zeta_U 1_U$. Let $\zeta$ be such a vector, and observe that for all $t > 0$, $F^{(\mathcal{V})}(\xi^{(\mathcal{V})} + t\zeta) - F^{(\mathcal{V})}(\xi^{(\mathcal{V})}) = F\left(\sum_{U \in \mathcal{V}} (\xi_U^{(\mathcal{V})} + t\zeta_U) 1_U\right) - F\left(\sum_{U \in \mathcal{V}} \xi_U^{(\mathcal{V})} 1_U\right) = F(x + td^{(x)}) - F(x)$. We deduce that $F^{(\mathcal{V})}$ admits a directional derivative at $\xi^{(\mathcal{V})}$ in direction $\zeta$ and that $F^{(\mathcal{V})'}(\xi^{(\mathcal{V})}, \zeta) = F'(x, d^{(x)})$. If $x$ is not a stationary point for $F$, then by proposition 2.4, $F'(x, d^{(x)}) < 0$, contradicting that $\xi^{(\mathcal{V})}$ is a stationary point for $F^{(\mathcal{V})}$.

Now $|\mathcal{V}| \leq |V|$, thus by contraposition the algorithm terminates correctly in a finite number of iterations. ∎

## 2.3    Implementation Considerations

As with the regularization of differentiable functionals, the reduced problem P2 presents the same structure as the main problem P1 in the presence of a separable nondifferentiable part; indeed, it decomposes as $\xi \mapsto \sum_{U \in \mathcal{V}} \sum_{v \in U} g_v(\xi_U) = \sum_{U \in \mathcal{V}} \gamma_U(\xi_U)$, where each $\gamma_U \colon \mathbb{R} \mapsto \,]-\infty, +\infty]$. Thus, any algorithm solving the reduced problems can also solve the main problem. Nevertheless, it often happens that the former run much faster and give more precise results than the latter, and the cut-pursuit can leverage this. When solutions with large constant connected components exist, performance is improved by orders of magnitude, as we show numerically in § 3. Interestingly, we observe on these problems that both crucial steps, namely solving the reduced problem and refining the partition, take a significant share of the total computational load. We give here practical implementation considerations, important for robustness and efficiency.

### 2.3.1    Discontinuities and Solutions of Reduced Problems

The cut-pursuit algorithm relies on directional derivatives of nondifferentiable functions, thus problems of discontinuity are to be expected. For once, the definition of the steepest ternary problem P4 at a point $x$ depends on the set of edges whose vertices share *exactly* the same value, $E_=^{(x)} \stackrel{\text{def}}{=} \{(u, v) \in E \mid x_u = x_v\}$; the absolute differences over all other edges are treated as differentiable terms at $x$, however small the difference is. Likewise, the computation of $\delta^+(x)$ and $\delta^-(x)$ requires identifying which coordinates $x_v$ lie at points of nondifferentiability of the corresponding $g_v$.

This is an important limitation, since in most applications, solutions of reduced problems are computed with iterative algorithms which are only asymptotically convergent; they are only approximate solutions, within a certain tolerance error. Such solutions sometimes even lie slightly outside the domain of the objective functionals, making things even worse. These facts cannot be ignored because nondifferentiability points are usually points of interest for the problem at hand; actually the very reason why nondifferentiable terms are considered in the first place. Consequently, we recommend setting up a threshold distance in coherence with the tolerance error specified for reduced problems. Coordinates which fall within this threshold of a nondifferentiability point are treated as if they were exactly at this differentiability point. When doing so, it is also profitable to merge together neighboring components in $\mathcal{V}$ which are assigned close values, because the lower the cardinal of $\mathcal{V}$, the smaller the reduced graph $\mathcal{G}$ and the faster the solutions of the reduced problems.

Observe that calling on such approximations, optimality considerations of § 2.2 are not strictly valid anymore; in particular, neighboring components can be alternatively merged after the reduced problem and split again after the steepest direction problem, *ad infinitum*. Thus, we also advocate terminating the algorithm when the iterate evolution is below a certain threshold, once again in coherence with the tolerance on the reduced problems.

Another numerical difficulty which is worth mentioning is that components of $\mathcal{V}$ can be very different in size, leading to bad conditioning of the reduced problem because large components have much more importance than small ones, preventing accurate estimation of the latter. A method that allows dealing with bad conditioning is thus required for solving the reduced problem.

In our numerical experiments, we use the preconditioning of the forward-Douglas–Rachford splitting algorithm illustrated by Raguet (2017), showing favorable behavior with respect to the above considerations on the problems that we consider.

### 2.3.2 Maximum Flow

Although different strategies have been developed for finding maximum flows in graphs, we only considered the augmenting path strategy of Boykov and Kolmogorov (2004), which seems well adapted to the structure of the flow graph $G_{\text{flow}}^{(x)}$ described in § 2.1 and figure 1. Let us underline that the horizontal structure of the flow graph is determined by the original graph $G$ and the components in the current partition $\mathcal{V}$; in particular, there is no horizontal edge between two different components of the partition, and a path from the source to the sink always goes through a unique component. This provides a natural way of parallelizing the computation of the maximum flow along the components. Moreover, the refinement of the partition which we propose is essentially hierarchical, each component being split into several parts, which also suggests a parallelization of the search for maximal connected components and might ease the memory structure. We have not implemented such parallelization so far, and leave it for future works.

In addition, the steepest ternary direction problem P4 can be solved by an alternative minimum cut strategy. For $d \in \mathbb{R}^V$, define respectively the coordinate-wise minimum and maximum $\min(d, 0), \max(d, 0) \in \mathbb{R}^V$, by for all $v \in V$, $\min(d, 0)_v \stackrel{\text{def}}{=} \min(d_v, 0)$ and $\max(d, 0)_v \stackrel{\text{def}}{=} \max(d_v, 0)$. Then, it can be shown from proposition 2.2 that for all $x \in \text{dom } F$,

$$F'(x, d) = F'(x, \min(d, 0)) + F'(x, \max(d, 0)) \,,$$

so that $d \mapsto F'(x, d)$ is minimized over $\{-1, 0, +1\}^V$ by the sum of a minimizer over $\{-1, 0\}^V$ and of a minimizer over $\{0, +1\}^V$. Each of the latter minima can be found by a minimum cut in an adapted flow graph like the one of figure 1, but with only one stage. This can be used to reduce memory requirements; alternatively, if memory is not a concern, the two minimizations could be performed in parallel. On our experiments below, we implemented this serially, with substantial gain in terms of memory and no loss in terms of running time.

## 2.4   Extension to Multidimensional Values

The very idea of the cut-pursuit algorithm 1 can be summarized as follows: solving a reduced problem on a partition of $V$, finding a steepest descent direction within a set $D$, and refining the partition accordingly. In theory, this strategy could be applied in any setting; however, if the non-differentiable part of $F$ besides the graph total variation is not a separable sum of unidimensional functionals, two difficulties arise. First, the set of descent directions $D$ necessary for obtaining an optimality certificate as in § 2.2 might be infinite. Second, even if $D$ is finite, the problem of finding the steepest descent direction might not be tractable.

Nevertheless, one can think of situations where these problems can be heuristically addressed. A typical one is when the nonsmooth functionals are not sums of unidimensional functionals, but are still separable over the graph $G$, in the sense that there is no edge between the coordinates over which each one is defined. This situation is better modeled by saying that the values at the vertices are multidimensional, say in $\mathbb{R}^K$ where $K$ is a finite set. The absolute value in the graph total variation can be replaced by any norm over $\mathbb{R}^K$, and the resulting objective functional is then defined, for all $x \in \mathbb{R}^{V \times K}$, as

$$F(x) \stackrel{\text{def}}{=} f(x) + \sum_{v \in V} g_v(x_v) + \sum_{(u,v) \in E} w_{(u,v)} \|x_u - x_v\|,$$

where now for all $v \in V$, $x_v \stackrel{\text{def}}{=} \left(x_{(v,k)}\right)_{k \in K} \in \mathbb{R}^K$ and $g_v \colon \mathbb{R}^K \to \,]{-}\infty, +\infty]$.

Compared to the setting of § 2.1, positive homogeneity of directional derivatives of the $g_v$ still holds but unit descent directions cannot be summarized by ascending, +1, or descending, −1: as soon as $|K| \geq 2$, there is an infinity of unit vectors. However, for a given vertex, only a handful of descent directions in $\mathbb{R}^K$ might seem relevant for the problem. Our first heuristic is to restrict the set of considered directions by choosing them greedily for each vertex. For example, if for each vertex $v$ only one direction $\bar{d}_v \in \mathbb{R}^K$ is considered, the set of directions is the Cartesian product $D \stackrel{\text{set}}{=} \times_{v \in V} \{0, \bar{d}_v\} \subset \mathbb{R}^{V \times K}$, and the corresponding steepest descent direction problem is binary. It is easy to show that for all $x \in \operatorname{dom} F$ and $d \in D$,

$$F'(x, d) = \sum_{\substack{v \in V \\ d_v = \bar{d}_v}} \delta(x, \bar{d}_v) + \sum_{(u,v) \in E_{=}^{(x)}} w_{(u,v)} \|d_u - d_v\|,$$

where the $\delta(x, \bar{d}_v)$ does not depend on $d$; so the problem can again be solved by finding a minimum cut in a (single stage) flow graph according to Theorem 4.1 of Kolmogorov and Zabih (2004), where condition (7) reduces to the triangle inequality for the norm defining the total variation.

It must be underlined here that the set $D$ above might be different at each iteration, depending on the current iterate $x$. Moreover, one can consider richer sets of direction per vertex, $D \stackrel{\text{set}}{=} \times_{v \in V} D_v$, where each $D_v$ is a finite subset of $\mathbb{R}^K$. Now, the steepest descent direction problem is a multilabel one, and in general cannot be easily solved. Fortunately, greedy strategies such as *α-expansion* or *α-β swap*, as described for instance by Boykov et al. (2001), can provide satisfactory approximate solutions by solving a succession of a few binary problems like the above.

We say that these approaches are heuristics because in the general case no optimality can be provided, neither for the original optimization problem, nor for the steepest descent problem when more than two descent directions are considered per vertex. Nonetheless, we show below, on a simplex-constrained labeling problem, that they can be efficient.

## 3 Numerical Experiments

Raguet (2017) illustrates his preconditioning of the forward-Douglas–Rachford splitting algorithm (PFDR) on medium- and large-scale problems arising respectively from signal processing and machine learning tasks, on which it compares favorably with state-of-the-art proximal splitting methods. We show the considerable improvement offered by the cut-pursuit (CP) approach on the exact same optimization problems, using PFDR for solving the reduced problems. In the comparisons, we also include the preconditioned primal-dual splitting algorithm of Pock and Chambolle (2011, PPD) because of its popularity; note that it is closely related to the alternating direction method of multipliers, often coined ADMM.

The experimental setting is extensively described by Raguet (2017, § 4), and we refer the reader to this note for details. We use slight improvements of his implementation of PFDR and PPD in C++ with parallelization of most operators with OpenMP specifications, and run the experiments on a personal computer with eight cores at 2.80 GHz and dual-channel DDR4 memory at 2.40 GHz. The source code for CP and PFDR is available at one of the author's GitHub repository.[1]

In the following, if $C$ is a convex closed set of a vector space $\Omega$, we note the *convex indicator* functional $\iota_C \colon \Omega \to\, ]-\infty, +\infty]\colon x \mapsto 0$ if $x \in C$, $+\infty$ otherwise.

### 3.1 Inverse Problem in Electroencephalography

Electroencephalography records brain activity via electrodes put at the surface of a subject's head. The relationship between activation of the brain regions and the electrodes' recording can be modeled by a linear operator called *lead-field* operator. The brain regions are modeled as vertices of a tridimensional mesh, $G \stackrel{\text{set}}{=} (V, E)$, and a brain activation map is thus a vector of $\mathbb{R}^V$. Yet, the number of electrodes being much smaller (here, $N \stackrel{\text{set}}{=} 91$) than the resolution of the desired brain image (here, $|V| \stackrel{\text{set}}{=} 19\,626$), the problem of retrieving brain activation map from the electrodes' recording is ill-posed. Moreover, the latter usually suffers from acquisition noise.

Fortunately, following Becker et al. (2014), a reasonable assumption is that at a given time, only scarce regions of the brain are really activated, and that spatially neighboring regions are often similarly activated. In addition, we use a recording time point where the entire signal is known to be nonnegative. All this prior knowledge can be enforced by modeling the brain source as a minimizer over $\mathbb{R}^V$ of

$$F \colon x \mapsto \tfrac{1}{2} \|y - \Phi x\|^2 + \sum_{v \in V} (\lambda_v |x_v| + \iota_{\mathbb{R}_+}(x_v)) + \sum_{(u,v) \in E} w_{(u,v)} |x_u - x_v|\,,$$

where $y \in \mathbb{R}^N$ is the observation over $N$ electrodes and $\Phi \colon \mathbb{R}^V \to \mathbb{R}^N$ is the lead-field operator. The first term is a square Euclidean norm ensuring coherence with the observation; it is differentiable. The second term is comprised of both a weighted $\ell_1$-norm and a convex indicator, enforcing respectively sparsity and positivity; it is nondifferentiable but separable over $G$. The third term is the graph total variation enforcing spatial similarity.

Altogether, this is of the form of problem P1, and the cut-pursuit algorithm can be easily applied following §§ 2.1 and 2.3. Once again, we refer the reader to the note of Raguet (2017, § 4) for details on the competing algorithms. Following his methodology, we prescribe stopping criteria as minimum relative evolution of the iterates, decreasing from $10^{-4}$ to $10^{-6}$; for the reduced

---

[1] `https://github.com/1a7r0ch3/CP_PFDR_graph_d1`

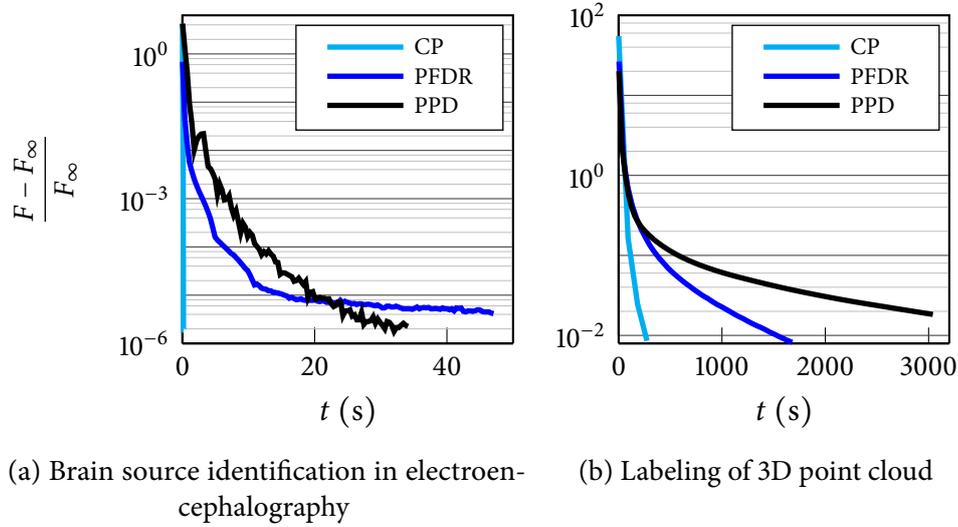

(a) Brain source identification in electroencephalography

(b) Labeling of 3D point cloud

Figure 4: Optimization comparisons.

Table 1: Brain source identification in electroencephalography: prediction performance and running time comparisons.[1] For the $10^{-8}$ stopping criterion, PPD and PFDR were stopped at $10^5$ iterations.

| stop. criterion | CP | | | PFDR | | | PPD | | |
|---|---|---|---|---|---|---|---|---|---|
| | DS | $DS_a$ | time (s) | DS | $DS_a$ | time (s) | DS | $DS_a$ | time (s) |
| $10^{-4}$ | 0.32 | 0.78 | 0.17 | 0.24 | 0.76 | 3 | 0.13 | 0.66 | 7 |
| $10^{-5}$ | 0.32 | 0.78 | 0.17 | 0.31 | 0.74 | 12 | 0.25 | 0.78 | 22 |
| $10^{-6}$ | 0.32 | 0.78 | 0.19 | 0.31 | 0.78 | 47 | 0.30 | 0.78 | 34 |
| $10^{-8}$ | 0.32 | 0.78 | 0.25 | 0.32 | 0.78 | 191 | 0.32 | 0.78 | 180 |

problems in CP, the stopping criterion is set to one thousandth of this value. We also consider longer runs of the algorithms with a stopping criterion of $10^{-8}$ for CP and stopped after $10^5$ iterations for PPD and PFDR.

In this experiment, and for all stopping criteria considered, CP terminates after 11 iterations, with only 20 maximal constant connected components. In such a favorable case, it outperforms the two other algorithms by several orders of magnitude, as illustrated on figure 4(a), where the optimal $F_\infty$ is approximated with CP with stopping criterion $10^{-8}$. For more in-depth comparison, we also report the computing times[1] in table 1. Moreover, as the data are synthetic, the original brain activity is known, and we can assess the relevance of the model for *brain source identification* by computing the *Dice score* between the supports of the retrieved activity and of the ground truth. We also report an approximate Dice score, $DS_a$, where small absolute values of the solutions are discarded with a simple 2-means algorithm.

---

[1] Differences with earlier results of Raguet (2017, § 4) are explained by better hardware and better memory management in the implementation, favorable to PFDR and even more to PPD.

## 3.2   Semantic Labeling of 3D Point Cloud

We consider the task of assigning a semantic label (car, vegetation, road, etc…) to each point of a 3D point cloud acquired with a LiDAR. This is usually performed with a supervised classifier such as a random forest, whose features can be derived from the local neighborhood of the points, or from the global structure of the scene; see for instance the works of Weinmann et al. (2015) and Guinard and Landrieu (2017).

If $V$ denotes the set of points and $K$ the set of labels, the random forest classifier provides a probabilistic classification $q \in \mathbb{R}^{V \times K}$, where for each $v \in V$, $q_v \stackrel{\text{def}}{=} (q_{(v,k)})_{k \in K}$ belongs to the simplex $\triangle_K \stackrel{\text{def}}{=} \{p \in \mathbb{R}^K \mid \sum_{k \in K} p_k = 1 \text{ and } \forall\, k \in K, p_k \geq 0\}$. Although it generally gives good results, it lacks the spatial regularity which can be expected from LiDAR acquisitions; following Landrieu et al. (2017), this can be improved by encoding an adjacency structure on a graph $G \stackrel{\text{set}}{=} (V, E)$, and minimizing over $\mathbb{R}^{V \times K}$ the functional

$$F \colon p \mapsto \sum_{v \in V} \mathrm{KL}(\beta u + (1-\beta) q_v, \beta u + (1-\beta) p_v) + \sum_{v \in V} \iota_{\triangle_K}(p_v) + \sum_{(u,v) \in E} w_{(u,v)} \|p_u - p_v\|_1 \,,$$

where for all $r, s \in \triangle_K$, $\mathrm{KL}(r, s) \stackrel{\text{def}}{=} \sum_{k \in K} r_k \log(r_k/s_k)$ is the *Kullback–Leibler divergence*, $u \stackrel{\text{def}}{=} (1/|K|)_{k \in K} \in \triangle_K$ is the uniform discrete distribution, and $\beta \in\, ]0, 1[$ is a small smoothing parameter. The first term favors similarity with the original predictions; it is differentiable. The second term ensures that each labeling is a discrete probability distribution; it is nondifferentiable but separable over $G$. The third term is the graph total variation enforcing spatial similarity, where we use the $\ell_1$ norm, well adapted to simplex-constrained values.

Altogether, this satisfies the multidimensional setting described in § 2.4. To construct the set of candidate descent directions, consider an iterate $p \in \triangle_K^V$ and a direction $d \in \mathbb{R}^{V \times K}$. Note that that for all $v \in V$, $\iota_{\triangle_K}'(p_v, d_v) = +\infty$ if $\sum_{k \in K} d_{(v,k)} \neq 0$, or if there exists $k \in K$ such that, either $p_{(v,k)} = 0$ and $d_{(v,k)} < 0$, or $p_{(v,k)} = 1$ and $d_{(v,k)} > 0$. Consequently, for each $v \in V$, we propose to define $k_v \in \arg\max_{k \in K} \{p_{(v,k)}\}$ a label with maximum probability, and set $D_v \stackrel{\text{set}}{=} \{0\} \cup \{\mathbf{1}_{\{k\}} - \mathbf{1}_{\{k_v\}} \in \mathbb{R}^K \mid k \in K \smallsetminus \{k_v\}\}$. The steepest descent direction is then a combinatorial problem with $|K|$ labels which we approximately solve with a single $\alpha$-expansion cycle.

The graph contains $|V| = 3\,000\,111$ vertices and $|E| = 17\,206\,938$, and the task comprises $|K| = 6$ classes. Stopping criteria are again taken from the experiments of Raguet (2017), and an estimate of the optimal value $F_\infty$ is computed with a longer run. Figure 4(b) represents the evolution of the objective functional values over time. The results are less impressive than in the previous experiment, but once again, CP reaches lower objective values an order of magnitude faster than PFDR; after only 4 iterations, with a total of 863 maximal constant connected components.

Let us underline that in this setting, the majority of the computational time is devoted to graph cuts. Indeed, starting at the direction $d = 0$, an entire $\alpha$-expansion cycle requires 5 successive graph cuts, over the huge original graph. There is thus room for significant improvements by parallelizing the cuts as explained along § 2.3.2, or by exploring better strategies for searching descent directions.

## 4   Conclusion and Perspectives

This paper provides a theoretical and practical framework for harnessing the speed of efficient graph-cut algorithms for a large class of graph-structured problems involving nondifferentiable terms alongside the total variation.

We believe that our algorithm overcomes three important limitations. First, solving total-variation regularized problems in high dimension is known to be difficult. Computational limitations might have led some works to use unconverged solutions, providing unsatisfying or inconsistent results. Cut-pursuit addresses this problem through its considerable acceleration, at least when the number of final constant connected components is reasonable.

Second, even when satisfying solutions can be found, practical applications often require lengthy exploration of regularization parameters at a prohibitive computational cost. Cut-pursuit can benefit from warm-restart of the partition, for scanning from high to low regularization strength, as already pointed out by Landrieu et al. (2017, § 2.6).

Third, convexity of the total-variation, while being convenient for optimization considerations, makes it sometimes not restrictive enough as a regularizer, admitting several solutions with many level sets. In some cases it is preferable to obtain spatially homogeneous solutions with only few level sets, which are in general better enforced with nonconvex regularizations. On the basis of its very principle and of our first numerical experiments, we argue that the cut-pursuit scheme favors the solutions with the fewest constant connected components, mitigating this third concern.

Altogether, it seems that many applications of the total-variation would benefit from our approach, and that our algorithm might spark a renewed interest of this regularization in the future.

# References

F. Bach, R. Jenatton, J. Mairal, and G. Obozinski. Optimization with sparsity-inducing penalties. *Foundations and Trends in Machine Learning*, 4(1):1–106, 2012.

H. Becker, L. Albera, P. Comon, R. Gribonval, and I. Merlet. Fast, variation-based methods for the analysis of extended brain sources. In *European Signal Processing Conference*, 2014.

Y. Boykov and V. Kolmogorov. An experimental comparison of min-cut/max-flow algorithms for energy minimization in vision. *IEEE Transactions on Pattern Analysis and Machine Intelligence*, 26(9):1124–1137, 2004.

Y. Boykov, O. Veksler, and R. Zabih. Fast approximate energy minimization via graph cuts. *IEEE Transactions on Pattern Analysis and Machine Intelligence*, 23(11):1222–1239, 2001.

A. Chambolle and J. Darbon. On total variation minimization and surface evolution using parametric maximum flows. *International Journal of Computer Vision*, 84(3):288–307, 2009.

P. L. Combettes and J.-C. Pesquet. A proximal decomposition method for solving convex variational inverse problems. *Inverse problems*, 24(6):65014–65040, 2008.

C. Couprie, L. Grady, L. Najman, J.-C. Pesquet, and H. Talbot. Dual constrained TV-based regularization on graphs. *SIAM Journal on Imaging Sciences*, 6(3):1246–1273, 2013.

S. Durand, J. Fadili, and M. Nikolova. Multiplicative noise removal using $\ell_1$ fidelity on frame coefficients. *Journal of Mathematical Imaging and Vision*, 36(3):201–226, 2010.

A. Gramfort, B. Thirion, and G. Varoquaux. Identifying predictive regions from fMRI with TV-$\ell_1$ prior. In *Pattern Recognition in Neuroimaging*. IEEE, 2013.

S. Guinard and L. Landrieu. Weakly supervised segmentation-aided classification of urban scenes from 3D LiDAR point clouds. *ISPRS Archives of the Photogrammetry, Remote Sensing and Spatial Information Sciences*, 2017.

Z. Harchaoui, A. Juditsky, and A. Nemirovski. Conditional gradient algorithms for norm-regularized smooth convex optimization. *Mathematical Programming*, 152(1-2):75–112, 2015.

J. Hiriart-Urruty and C. Lemaréchal. *Fundamentals of Convex Analysis*. Grundlehren Text Editions. Springer Berlin Heidelberg, 2004.

H. Ishikawa. Exact optimization for Markov random fields with convex priors. *IEEE transactions on pattern analysis and machine intelligence*, 25(10):1333–1336, 2003.

V. Kolmogorov and R. Zabih. What energy functions can be minimized via graph cuts? *IEEE Transactions on Pattern Analysis and Machine Intelligence*, 26(2):147–159, 2004.

L. Landrieu and G. Obozinski. Cut pursuit: Fast algorithms to learn piecewise constant functions on general weighted graphs. *SIAM Journal on Imaging Sciences*, 10(4):1724–1766, 2017.

L. Landrieu, H. Raguet, B. Vallet, C. Mallet, and M. Weinmann. A structured regularization framework for spatially smoothing semantic labelings of 3D point clouds. *Journal of Photogrammetry and Remote Sensing*, 132:102–118, 2017.

C. Nieuwenhuis, E. Töppe, and D. Cremers. A survey and comparison of discrete and continuous multi-label optimization approaches for the Potts model. *International journal of computer vision*, 104(3):223–240, 2013.

M. Nikolova. A variational approach to remove outliers and impulse noise. *Journal of Mathematical Imaging and Vision*, 20(1-2):99–120, 2004.

N. Omranian, J. M. Eloundou-Mbebi, B. Mueller-Roeber, and Z. Nikoloski. Gene regulatory network inference using fused LASSO on multiple data sets. *Scientific reports*, 6, 2016.

T. Pock and A. Chambolle. Diagonal preconditioning for first order primal-dual algorithms in convex optimization. In *IEEE International Conference on Computer Vision*, pages 1762–1769. IEEE, 2011.

H. Raguet. A note on the forward-Douglas–Rachford splitting for monotone inclusion and convex optimization. *preprint*, 2017.

H. Raguet and L. Landrieu. Preconditioning of a generalized forward-backward splitting and application to optimization on graphs. *SIAM Journal on Imaging Sciences*, 8(4):2706–2739, 2015.

T. Takayama and A. Iwasaki. Optimal wavelength selection on hyperspectral data with fused LASSO for biomass estimation of tropical rain forest. *ISPRS Annals of Photogrammetry, Remote Sensing and Spatial Information Sciences*, pages 101–108, 2016.

R. Tibshirani, M. Saunders, S. Rosset, J. Zhu, and K. Knight. Sparsity and smoothness via the fused lasso. *Journal of the Royal Statistical Society: Series B (Statistical Methodology)*, 67(1): 91–108, 2005.

M. Weinmann, B. Jutzi, S. Hinz, and C. Mallet. Semantic point cloud interpretation based on optimal neighborhoods, relevant features and efficient classifiers. *ISPRS Journal of Photogrammetry and Remote Sensing*, 105:286–304, 2015.

X. Wu, J. Zheng, Y. Cai, and C.-W. Fu. Mesh denoising using extended ROF model with $\ell_1$ fidelity. In *Computer Graphics Forum*, volume 34, pages 35–45. Wiley Online Library, 2015.

B. Xin, Y. Kawahara, Y. Wang, L. Hu, and W. Gao. Efficient generalized fused LASSO and its applications. *ACM Transactions on Intelligent Systems and Technology*, 7(4):1–22, 2016.